\documentclass[11pt,twoside]{article}

\usepackage{a4wide}
\usepackage{amsfonts}
\usepackage{amssymb}
\usepackage{amsmath}
\usepackage{graphicx}

\newcommand{\ignore}[1]{}

\def\@begintheorem#1#2{\par\bgroup{\sc #1\ #2. }\it\ignorespaces}
\def\@opargbegintheorem#1#2#3{\par\bgroup{\sc #1\ #2\ (#3). } \it\ignorespaces}
\def\@endtheorem{\egroup}
\newtheorem{theorem}{Theorem}[section]
\newtheorem{corollary}[theorem]{Corollary}
\newtheorem{lemma}[theorem]{Lemma}

\newtheorem{example}[theorem]{Example}
\newtheorem{proposition}[theorem]{Proposition}
\newtheorem{definition}[theorem]{Definition}
\newcommand{\bt}[1]{\begin{theorem}\label{#1}}
\newcommand{\bc}[1]{\begin{corollary}\label{#1}}
\newcommand{\bl}[1]{\begin{lemma}\label{#1}}
\newcommand{\be}[1]{\begin{example}\label{#1}}
\newcommand{\bp}[1]{\begin{proposition}\label{#1}}
\newcommand{\ba}[1]{\begin{algorithm}\rm\label{#1}}
\newcommand{\bd}[1]{\begin{definition}\rm\label{#1}}{\normalsize }
\newcommand{\bpr}{\noindent {\em Proof. }}
\newcommand{\et}{\end{theorem}}
\newcommand{\ec}{\end{corollary}}
\newcommand{\el}{\end{lemma}}
\newcommand{\ee}{\end{example}}
\newcommand{\ep}{\end{proposition}}
\newcommand{\ed}{\end{definition}}
\newcommand{\epr}{{\ \vbox{\hrule\hbox{%
\vrule height1.3ex\hskip0.8ex\vrule}\hrule}}\\\par}

\def\R{\mathbb{R}}
\def\Z{\mathbb{Z}}
\def\B{{\cal B}}
\def\D{{\cal D}}
\def\E{{\cal E}}
\def\conv{\mbox{conv}}

\begin{document}

\title{\bf Degree Sequence Optimization\\ and\\ Extremal Degree Enumerators}

\author{Shmuel Onn}

\date{}

\maketitle

\begin{abstract}
The degree sequence optimization problem is to find a subgraph of a given graph which
maximizes the sum of given functions evaluated at the subgraph degrees.
Here we study this problem by replacing degree sequences, via suitable nonlinear transformations,
by suitable degree enumerators, and we introduce suitable degree enumerator polytopes.

We characterize their vertices, that is, the extremal degree enumerators, for complete
graphs and some complete bipartite graphs, and use these characterizations to obtain
simpler and faster algorithms for optimization over degree sequences for such graphs.

\vskip.2cm
\noindent {\bf Keywords:} graph, degree sequence, combinatorial optimization, polytope
\end{abstract}

\section{Introduction}

The {\em degree sequence} of a simple $n$-graph $G=(V,E)$ with $V=\{v_1,\dots,v_n\}$
is the vector $d(G)=(d_{v_1}(G),\dots,d_{v_n}(G))$, where $d_{v_i}(G):=|\{e\in E:v_i\in e\}|$
is the degree of vertex $v_i$ for all $i$. Degree sequences have been studied
by many authors, starting from their celebrated effective characterization by
Erd\H{o}s and Gallai \cite{EG}, see e.g. \cite{EKM} and the references therein.
Throughout, by a {\em subgraph} $G\subseteq H$ of a given graph $H$ on $V$ we mean a subgraph
of $H$ with the same vertex set $V$. We also use throughout the somewhat nonstandard
notation $[n]:=\{0,1,\dots,n-1\}$. We continue the study of optimization problems over
degree sequences, e.g. \cite{AS,DLMO,KO,PPS} and references therein, of finding a subgraph
of a given graph maximizing the sum of given functions evaluated at the subgraph degrees.
Previous related work also considered polytopes of degree sequences, e.g. \cite{Kor,MS,PS}.
Here we replace degree sequences, via suitable nonlinear transformations, by suitable
{\em degree enumerators}, and introduce suitable {\em degree enumerator polytopes}.
We characterize their vertices for complete graphs and some complete bipartite graphs,
and use these characterizations to obtain faster algorithms for optimization over degree
sequences for such graphs. Following is a more detailed description.

\subsection{Graphs}

First, we consider the following problem, with one function evaluated on all vertex degrees.

\vskip.2cm\noindent{\bf Degree Sequence Optimization.}
Given an $n$-graph $H$ with vertex set $V$, and a function $f:\{0,\dots,n-1\}\rightarrow\Z$,
find a subgraph $G\subseteq H$ which attains maximum $\sum_{i=1}^n f(d_{v_i}(G))$.

\vskip.2cm
It extends the perfect matching problem: indeed, $H$ has a perfect matching if and only if the
optimal value of the degree sequence optimization problem over $H$ and $f(x):=-(x-1)^2$ is zero.
In fact, the problem is generally NP-hard. To see this, for a given $n$, define
$$f(x)\ :=\ \left\{\begin{array}{ll}
n-1, & x=0; \\
n, & x=3; \\
0, & \hbox{otherwise.}
\end{array}\right.\ .$$
If a given $H$ has a nonempty cubic subgraph $G$ (that is, with every vertex of degree $0$ or $3$)
then $\sum_{i=1}^n f(d_{v_i}(G))\geq n^2-n+1$, whereas if not, then any $G\subseteq H$
satisfies $\sum_{i=1}^n f(d_{v_i}(G))\leq n^2-n$. So optimization enables to decide
nonempty cubic subgraphs, which is NP-complete \cite{GJ}.

\vskip.2cm
In Section 2 we consider an equivalent formulation of the problem via a nonlinear
transformation as follows. We let the {\em degree enumerator} of $G$ be the
vector $e(G)=(e_0(G),\dots,e_{n-1}(G))$, where
$e_i(G):=|\{v_j\in V:d_{v_j}(G)=i\}|$ is the number of vertices of $G$ of degree $i$.
For instance, the $7$-graph $G$ with $V=\{1,\dots,7\}$ and $E=\{12,16,17,23,34,37,45,56,57\}$
has degree sequence $d(G)=3232323$ and degree enumerator $e(G)=0034000$.
We index vectors in $\R^n$ as $x=(x_0,\dots,x_{n-1})$ which is consistent with the
indexation of degree enumerators. In particular, the unit vectors are denoted by
${\bf 1_i}$ for $i=0,1,\dots,n-1$. We identify functions $f:\{0,\dots,n-1\}\rightarrow\Z$
with vectors $f=(f_0,\dots,f_{n-1})$ via $f_i:=f(i)$. The standard inner product of
two vectors $f,z\in\R^n$ is denoted $fz:=\sum_{i=0}^{n-1}f_iz_i$.
Thus, given graph $G$ and function $f$ with corresponding vector $f$,
the objective function value of $G$ under $f$ becomes
$$\sum_{j=1}^nf(d_{v_j}(G))\ =\ \sum_{i=0}^{n-1}f_ie_i(G)\ =\ fe(G)\ .$$

With this, we have the following equivalent formulation of degree sequence optimization.

\vskip.2cm\noindent{\bf Degree Sequence Optimization (equivalent formulation).}
Given an $n$-graph $H$ and a vector $f\in\Z^n$, find a subgraph $G\subseteq H$
which attains maximum inner product value $fe(G)$.

\vskip.4cm
We study the following problem, of characterizing the {\em extremal degree enumerators}, which
in particular has implications to faster algorithms for degree sequence optimization.

\vskip.2cm\noindent{\bf Extremal Degree Enumerator Problem.}
Given an $n$-graph $H$, characterize the degree enumerators $e\in\Z^n$ for which
there is a vector $f\in\Z^n$, such that every subgraph $G\subseteq H$ optimal to the
degree optimization problem over $H$ and $f$, has degree enumerator $e(G)=e$.

\vskip.2cm
To study this problem we introduce the following natural polytope associated with $H$.

\vskip.2cm\noindent{\bf Degree Enumerator Polytope.}
The {\em degree enumerator polytope} of an $n$-graph $H$ is
$$\E_H\ :=\ \conv\{e(G)=(e_0(G),\dots,e_{n-1}(G))\ :\ G\subseteq H\}\ \subset\ \R^n\ .$$

Then the extremal degree enumerators are precisely the vertices of $\E_H$ and the problem is to
characterize them. That is, we want to characterize the degree enumerators $e\in\Z^n$ for which
there is some $f\in\Z^n$ such that $e$ is the unique maximizer of $fz$ over the polytope $\E_H$.

\vskip.2cm
As mentioned, degree sequence optimization is NP-hard for general graphs, and therefore the
characterization of vertices of $\E_H$ for an arbitrary given graph $H$ is presumably
impossible. So it is interesting to study this problem for special classes of graphs $H$.
The first natural case is the complete graph $H:=K_n$ and we denote the corresponding polytope
by $\E_n:=\E_{K_n}$. It was shown in \cite[Theorem 1.2]{DLMO} that the degree sequence
optimization problem over $K_n$ can be solved by dynamic programming in polynomial $O(n^7)$ time.
Here, in Section 2, we solve the extremal degree enumerator problem over complete graphs
and completely characterize the vertices of $\E_n$ in Theorem 2.1,
each vertex $e$ along with a suitable graph $G$ with $e=e(G)$.\break
For instance, for $n=7$, brute force computation shows that there are $342$ degree enumerators,
but by Theorem 2.1 only $16$ out of them are vertices of $\E_7$ and are the columns of the matrix
$$\left(\begin{array}{cccccccccccccccc}
7 & 0 & 0 & 0 & 1 & 0 & 0 & 0 & 1 & 0 & 0 & 0 & 1 & 0 & 0 & 0\\
0 & 0 & 0 & 0 & 6 & 6 & 6 & 6 & 0 & 0 & 0 & 0 & 0 & 0 & 0 & 0\\
0 & 7 & 0 & 0 & 0 & 1 & 0 & 0 & 0 & 1 & 0 & 0 & 0 & 1 & 0 & 0\\
0 & 0 & 0 & 0 & 0 & 0 & 0 & 0 & 6 & 6 & 6 & 6 & 0 & 0 & 0 & 0\\
0 & 0 & 7 & 0 & 0 & 0 & 1 & 0 & 0 & 0 & 1 & 0 & 0 & 0 & 1 & 0\\
0 & 0 & 0 & 0 & 0 & 0 & 0 & 0 & 0 & 0 & 0 & 0 & 6 & 6 & 6 & 6\\
0 & 0 & 0 & 7 & 0 & 0 & 0 & 1 & 0 & 0 & 0 & 1 & 0 & 0 & 0 & 1
\end{array}\right)\ .
$$

The characterization of Theorem 2.1 leads in Corollary 2.2 to a much simpler algorithm for
doing degree sequence optimization over complete graphs, in much faster $O(n^2)$ time.

\subsection{Bipartite graphs}

Next we consider refined, bi-analogs, of all the above problems and notions for bipartite graphs,
as follows. We consider bipartite graphs $G=(U,V,E)$ with $U=\{u_1,\dots,u_m\}$
the set of {\em left} vertices, $V=\{v_1,\dots,v_n\}$ the set of {\em right} vertices,
and $E\subseteq U\times V$ the set of edges.
We consider the following problem, with two functions, one on the left and one on the right.

\vskip.2cm\noindent{\bf Degree Sequence Bi-Optimization.}
Given bipartite $(m,n)$-graph $H=(U,V,E)$ and functions $f:\{0,\dots,n\}\rightarrow\Z$,
$g:\{0,\dots,m\}\rightarrow\Z$, find a subgraph $G\subseteq H$ attaining maximum
$$\sum_{i=1}^m f(d_{u_i}(G))+\sum_{j=1}^n g(d_{v_j}(G))\ .$$

\vskip.2cm
This problem is again generally NP-hard. To see this, define the functions
$f(x):=x(x-3)$ and $g(x):=-(x-1)^2$. Given $3$-subsets $S_1,\dots,S_m$ of $V:=\{v_1,\dots,v_n\}$,
define a bipartite graph $H=(U,V,E)$ with $U:=\{u_1,\dots,u_m\}$ and $E:=\{u_iv_j:v_j\in S_i\}$.
Then the optimal value of the above bi-optimization problem over $H$, $f$, and $g$,
is zero if and only if some of the $S_i$ can be selected so as to form a partition of $V$,
which is NP-complete to decide \cite{GJ}.

\vskip.2cm
In Section 3, we consider an equivalent formulation of the above bi-problem, via the following
nonlinear transformation. We let the {\em degree bi-enumerator} of bipartite $(m,n)$-graph
$G=(U,V,E)$ be the vector $b(G)=a(G)\oplus c(G)$ in the direct sum $\R^{n+1}\oplus\R^{m+1}$, defined by
\begin{eqnarray*}
a_k(G)&:=&|\{i\in U:d_{u_i}(G)=k\}|,\quad k=0,\dots,n\ , \\
c_k(G)&:=&|\{j\in V:d_{v_j}(G)=k\}|,\quad k=0,\dots,m\ .
\end{eqnarray*}
For instance, the bipartite $(2,6)$-graph $G$ with $U=\{u_1,u_2\}$, $V=\{v_1,\dots,v_6\}$ and with
$E=\{u_1v_1,u_1v_2,u_2v_2,u_2v_3,u_2v_4,u_2v_5,u_2v_6\}$
has bi-enumerator $a(G)\oplus c(G)=0010010\oplus051$.
We denote vectors in the direct sum space $\R^{n+1}\oplus\R^{m+1}$ in the following equivalent ways,
\begin{eqnarray*}
(x,y)&=& x \oplus y\ =\ (x_0,\dots,x_n,y_0,\dots,y_m)\ =\ (x_0,\dots,x_n)\oplus(y_0,\dots,y_m) \\
&=& (x_0\cdot{\bf 1}_0+\cdots+x_n\cdot{\bf 1}_n)\oplus(y_0\cdot{\bf 1}_0+\dots+y_m\cdot{\bf 1}_m) \\
&=& (x_0\cdot{\bf 1}_0+\cdots+x_n\cdot{\bf 1}_n)\oplus (y_0,\dots,y_m)
\ =\ (x_0,\dots,x_n)\oplus(y_0\cdot{\bf 1}_0+\dots+y_m\cdot{\bf 1}_m)\ .
\end{eqnarray*}
Again we identify functions $f:\{0,\dots,n\}\rightarrow\Z$ and $g:\{0,\dots,m\}\rightarrow\Z$
with vectors $f=(f_0,\dots,f_n)$ via $f_i:=f(i)$ and  $g=(g_0,\dots,g_m)$ via $g_j:=g(j)$.
The standard inner product of two vectors $f\oplus g, a\oplus c\in\R^{n+1}\oplus\R^{m+1}$
is $\langle f\oplus g,a\oplus c\rangle=fa+gc=\sum_{i=0}^nf_ia_i+\sum_{j=0}^mg_jc_j$.
Thus, given bipartite $(m,n)$-graph $G$ with bi-enumerator $b(G)=a(G)\oplus c(G)$,
and functions $f,g$ with corresponding vectors $f,g$,
the objective function value of $G$ under $f$ and $g$ becomes
$$\sum_{i=1}^m f(d_{u_i}(G))+\sum_{j=1}^n g(d_{v_j}(G))\ =\ fa(G)+gc(G)
\ =\ \langle f\oplus g,a(G)\oplus c(G)\rangle\ =\ \langle f\oplus g,b(G)\rangle\ .$$

With this, we have the following equivalent formulation of degree sequence bi-optimization.

\vskip.2cm\noindent{\bf Degree Sequence Bi-Optimization (equivalent formulation).}
Given bipartite $(m,n)$-graph $H$ and $f\in\Z^{n+1}$, $g\in\Z^{m+1}$, find $G\subseteq H$
attaining maximum inner product $\langle f\oplus g,b(G)\rangle$.

\vskip.4cm
We study the following problem, of characterizing the {\em extremal degree bi-enumerators}, which
in particular has implications to faster algorithms for degree sequence bi-optimization.

\vskip.2cm\noindent{\bf Extremal Degree Bi-Enumerator Problem.}
Given bipartite $(m,n)$-graph $H$, characterize the degree bi-enumerators
$b=a\oplus c\in\Z^{n+1}\oplus\Z^{m+1}$ for which there are vectors $f\in\Z^{n+1}$ and $g\in\Z^{m+1}$,
such that every subgraph $G\subseteq H$ which is optimal to the degree sequence bi-optimization
problem over $H$, $f$, $g$, has degree bi-enumerator $b(G)=a(G)\oplus c(G)=a\oplus c=b$.

\vskip.2cm
To study this problem we introduce the following natural polytope associated with $H$.

\vskip.2cm\noindent{\bf Degree Bi-Enumerator Polytope.}
The {\em degree bi-enumerator polytope} of a bipartite $H$ is
$$\B_H\ :=\ \conv\{b(G)=a(G)\oplus c(G)\ :\ G\subseteq H\}
\ \subset\ \R^{n+1}\oplus\R^{m+1}\ .$$

Then the extremal degree bi-enumerators are precisely the vertices of $\B_H$ and the problem is to
characterize them. So we want to characterize bi-enumerators $b\in\Z^{n+1}\oplus\Z^{m+1}$
for which there are $f\in\Z^{n+1}$ and $g\in\Z^{m+1}$ such that $b$ is the unique maximizer of
$\langle f\oplus g,z\rangle$ over $\B_H$.

\vskip.2cm
As mentioned, degree sequence bi-optimization is NP-hard for general bipartite graphs,
so characterization of vertices of $\B_H$ for an arbitrary given bipartite $H$ is presumably
impossible. So it is interesting to study this problem for special classes of $H$.
The first natural case is the complete bipartite graph $H:=K_{m,n}$ and we denote the
corresponding polytope by $\B_{m,n}:=\B_{K_{m,n}}$. It was shown in \cite[Theorem 1.1]{KO}
that the degree sequence bi-optimization problem over $K_{m,n}$ can be solved by dynamic
programming in polynomial $O(m^5n^4)$ time. While $\B_{1,n}$ is just a simplex,
it turns out that already for $m=2$ the vertices of the polytope $\B_{2,n}$ have a much more
complicated structure than those of the enumerator polytope $\E_n$. Nonetheless, in Section 3
we settle the $m=2$ case and completely characterize the vertices of $\B_{2,n}$ for arbitrary $n$
in Theorem 3.3, each vertex $b$ along with a suitable bipartite graph $G$ with $b=b(G)$.\break
For instance, for $n=6$, there are $50$ degree bi-enumerators by Proposition
\ref{m=2_bi-enumerators}, but by Theorem 3.3 only $21$ out of them are vertices of $\B_{2,6}$
and are the columns of the matrix
$$
\left(\begin{array}{ccccccccccccccccccccc}
2 & 0 & 0 & 0 & 0 & 0 & 0 & 0 & 0 & 0 & 0 & 0 & 1 & 1 & 1 & 0 & 0 & 0 & 0 & 0 & 0\\
0 & 2 & 2 & 0 & 0 & 0 & 0 & 0 & 0 & 0 & 0 & 0 & 0 & 0 & 0 & 1 & 1 & 1 & 0 & 0 & 0\\
0 & 0 & 0 & 2 & 2 & 0 & 0 & 0 & 0 & 0 & 0 & 0 & 0 & 0 & 0 & 0 & 0 & 0 & 1 & 1 & 1\\
0 & 0 & 0 & 0 & 0 & 2 & 2 & 0 & 0 & 0 & 0 & 0 & 0 & 0 & 0 & 0 & 0 & 0 & 0 & 0 & 0\\
0 & 0 & 0 & 0 & 0 & 0 & 0 & 2 & 2 & 0 & 0 & 0 & 1 & 0 & 0 & 1 & 0 & 0 & 1 & 0 & 0\\
0 & 0 & 0 & 0 & 0 & 0 & 0 & 0 & 0 & 2 & 2 & 0 & 0 & 1 & 0 & 0 & 1 & 0 & 0 & 1 & 0\\
0 & 0 & 0 & 0 & 0 & 0 & 0 & 0 & 0 & 0 & 0 & 2 & 0 & 0 & 1 & 0 & 0 & 1 & 0 & 0 & 1\\
\oplus & \oplus & \oplus & \oplus & \oplus & \oplus & \oplus & \oplus & \oplus & \oplus &
\oplus & \oplus & \oplus & \oplus & \oplus & \oplus & \oplus & \oplus & \oplus & \oplus & \oplus \\
6 & 5 & 4 & 4 & 2 & 3 & 0 & 2 & 0 & 1 & 0 & 0 & 2 & 1 & 0 & 1 & 0 & 0 & 0 & 0 & 0\\
0 & 0 & 2 & 0 & 4 & 0 & 6 & 0 & 4 & 0 & 2 & 0 & 4 & 5 & 6 & 5 & 6 & 5 & 6 & 5 & 4\\
0 & 1 & 0 & 2 & 0 & 3 & 0 & 4 & 2 & 5 & 4 & 6 & 0 & 0 & 0 & 0 & 0 & 1 & 0 & 1 & 2
\end{array}\right)\ .
$$

The characterization of Theorem 3.3 leads in Corollary 3.4 to a much simpler algorithm
for doing degree sequence bi-optimization over $K_{2,n}$, in much faster $O(n^2)$ time.

\subsection{Overview}

Before proceeding, we remark on some related problems and results. First, the degree sequence
optimization and bi-optimization problems can be simultaneously extended as follows. Given now
is an $n$-graph $H$ on $V$ and $n$ functions $f^1,\dots,f^n:\{0,1,\dots,n-1\}\rightarrow\Z$.
The extended problem is to find $G\subseteq H$ attaining maximum
$\sum_{i=1}^n f^i(d_{v_i}(G))$. The complexity of this problem is intriguingly open
already for the complete graph $K_n$ and the complete bipartite graph $K_{m,n}$.
But partial results are known in various cases, most notably, the problem can be solved
in polynomial time over any graph $H$ when all functions are {\em concave} \cite{AS,DO}.
This problem also extends the general factor problem \cite{Cor}, where instead of functions $f^i$,
given are sets $D_i\subseteq[n]$, and we are to decide if there is a subgraph
$G\subseteq H$ with $d_{v_i}(G)\in D_i$ for all $i$.

A variant of the problem is as follows. For a degree sequence $d(G)=(d_{v_1}(G),\dots,d_{v_n}(G))$
let $\bar{d}(G)=(\bar{d}(G)_1,\dots,\bar{d}(G)_n)$ be the permutation of $d(G)$
satisfying $\bar{d}(G)_1\leq\cdots\leq\bar{d}(G)_n$. The problem is to find
$G\subseteq H$ maximizing $\sum_{i=1}^n f^i(\bar{d}(G)_i)$. This problem and a slight further
variant of it were shown to be polynomial time solvable by dynamic programming for complete graphs $K_n$
\cite[Theorem 1.2]{DLMO} and complete bipartite graphs $K_{m,n}$ \cite[Theorem 1.2]{KO},
respectively,
implying the aforementioned polynomial time solvability of degree optimization with a
single function $f$ over $K_n$ and bi-optimization with two functions $f,g$ over $K_{m,n}$.
We remark that in \cite[Theorem 1.2]{DLMO} only this result, for
identical functions, was stated, but it is clear that the proof also provides an
algorithm for $n$ functions applied to the permuted sequence $\bar{d}(G)$.

Another polytope which can be naturally associated with each $n$-graph $H$ on $V$
is its {\em degree sequence polytope} $\D_H:=\conv\{d(G):G\subseteq H\}$.
Its vertices correspond to degree sequences maximizing $\sum_{i=1}^nf_id_{v_i}(G)$
for some integers $f_1,\dots,f_n$, which is the very special case of the multi-function
objective $\sum_{i=1}^nf^i(d_{v_i}(G))$ where each $f^i$ is a {\em linear} function
$f^i(x):=f_ix$. The polytope $\D_n:=\D_{K_n}$ of the complete graph $H=K_n$ was studied
quite extensively in the literature. In particular, its vertices where characterized
in \cite{Kor} to be precisely the degree sequences of {\em threshold graphs} \cite{MP}.
See also \cite{MS,PS} for further work on this line of research.

Finally, we remark that all the above optimization problems and polytopes can be extended to
$k$-uniform hypergraphs for any $k$, but are harder to study, see for instance \cite{DLMO,MS}.

\section{Extremal enumerators of complete graphs}

Throughout this section we assume that the vertex set of an $n$-graph is $[n]=\{0,1,\dots,n-1\}$.
As usual, for positive integers $p,q$ we denote by $q\!\!\mod p$ the remainder of $q$ upon division
by $p$, for instance $14\!\!\mod 7=0$. For suitable values of $n$ and $r$ we define an
$r$-regular $n$-graph $G_n(r):=([n],E_n(r))$ as follows.
For any positive integer $n$ and any even $0\leq r\leq n-1$, let
$$E_n(r)\ :=\ \left\{\{i,(i+j)\!\!\!\!\mod n\}\ :\
i=0,\dots,n-1,\ \ j=1,\dots,{r\over 2}\right\}\ .$$
Arranging vertices on a circle, each vertex is connected to
the ${r\over 2}$ closest vertices on each side.
For instance, for $n=7$, $r=4$, we get $E_7(4)=\{01,02,12,13,23,24,34,35,45,46,56,50,60,61\}$.

For any even $n$ and any odd $1\leq r\leq n-1$,
we define $E_n(r)$ using the perfect matching
$$M_n\ :=\ \left\{\left\{i,{n\over 2}+i\right\}\ :\
i=0,\dots,{n\over 2}-1\right\},\quad\mbox{by}\quad E_n(r)\ :=\ E_n(r-1)\ \uplus\ M_n\ .$$
With vertices arranged on a circle again,
each vertex is matched in $M_n$ to the opposite one.
For instance, for $n=8$ and $r=3$, we get $E_8(3)=\{01,12,23,34,45,56,67,70,04,15,26,37\}$.

However, for any odd $n$ and any odd $r$, an $r$-regular $n$-graph does not exist, since the
sum of vertex degrees, which is always even, would have been $nr$, which is odd, so impossible.

Instead, for any odd $n$, any odd $0<r<n-1$, and any even $0\leq s\leq n-1$,
we define $G_n(r,s):=([n],E_n(r,s))$ as follows.
If $s=0$ then let $E_n(r,s):=E_{n-1}(r)$. If $s\geq 2$ then let
\begin{eqnarray*}
E_n(r,s) & := & E_{n-1}(r)
\ \setminus\ \left\{\left\{i,{n-1\over 2}+i\right\}:i=0,\dots,{s\over 2}-1\right\} \\
& \uplus    & \left\{\{i,n-1\}:i=0,\dots,{s\over 2}-1\right\}
\ \uplus\ \left\{\left\{{n-1\over 2}+i,n-1\right\}:i=0,\dots,{s\over 2}-1\right\}\ .
\end{eqnarray*}
Thus, $G_n(r,s)$ is obtained from $G_{n-1}(r)$ by removing $s\over 2$ of the edges
in the matching $M_{n-1}$ and connecting the end points of the removed matching edges
to a new vertex $n-1$. So $G_n(r,s)$ is almost $r$-regular,
in that vertex $n-1$ has degree $s$ and all other vertices have degree $r$.\break
For instance, for $n=7$, $r=3$, $s=4$, we get $E_7(3,4)=\{01,12,23,34,45,50,25,06,36,16,46\}$.

\vskip.2cm
We now show that the extremal degree enumerators of $K_n$ are exactly those of suitable $r$-regular
graphs $G_n(r)$ and almost $r$-regular graphs $G_n(r,s)$ according to the parities of $n$, $r$.

\bt{complete}
The vertices of the degree enumerator polytope $\E_n=\E_{K_n}$ are as follows.
\begin{itemize}
\item
For even $n$, they are $e_n(r):=n\cdot{\bf 1}_r=e(G_n(r))$ for $0\leq r\leq n-1$.
\item
For odd $n$, they are:
\begin{enumerate}
\item
$e_n(r)=n\cdot{\bf 1}_r=e(G_n(r))$ for even $0\leq r\leq n-1$;
\item
$e_n(r,s):=(n-1)\cdot{\bf 1}_r+{\bf 1}_s=e(G_n(r,s))$ for odd $0<r<n-1$ and even $0\leq s\leq n-1$.
\end{enumerate}
\end{itemize}
The number of vertices for even $n$ is $n$ and for odd $n$
is $({n+1\over2})^2$, so grows quadratically in $n$.
\et

\bpr
In one direction we show that any extremal degree enumerator, namely, any vertex of the polytope
$\E_n$, is one of the vectors indicated by the theorem. Let $v$ be any vertex of $\E_n$
and let $f\in\R^n$ be such that $v$ is the unique maximizer of $fz$ over $\E_n$.
Pick $0\leq r\leq n-1$ with $f_r=\max\{f_0,\dots,f_{n-1}\}$ and
pick $0\leq s\leq n-1$ even with $f_s=\max\{f_i:0\leq i\leq n-1,\ i\ \mbox{even}\}$.

First suppose $n$ is even. Let $x:=e(G_n(r))=n\cdot{\bf 1}_r$ be the degree enumerator of
the $r$-regular graph $G_n(r)$. Then for any degree enumerator $y=e(G)$ of any $n$-graph $G$ we have
$$fy\ =\ \sum_{i=0}^{n-1}f_iy_i\ \leq\ \sum_{i=0}^{n-1}f_ry_i
\ =\ f_r\sum_{i=0}^{n-1}y_i\ =\ nf_r\ =\ fx\ .$$
So $x$ maximizes $f$ over $\E_n$ and since $v$ is the
unique maximizer we have $v=x=n\cdot{\bf 1}_r$.

Next suppose $n$ is odd and suppose $f_s=f_r$.
Let $x:=e(G_n(s))=n\cdot{\bf 1}_s$ be the degree enumerator of $G_n(s)$.
Then for any degree enumerator $y=e(G)$ of any $n$-graph $G$ we have
$$fy\ =\ \sum_{i=0}^{n-1}f_iy_i\ \leq\ \sum_{i=0}^{n-1}f_ry_i
\ =\ f_r\sum_{i=0}^{n-1}y_i\ =\ nf_r\ =\ nf_s\ =\ fx\ .$$
So $x$ maximizes $f$ over $\E_n$ and since $v$ is the
unique maximizer we have $v=x=n\cdot{\bf 1}_s$.

Finally suppose $n$ is odd and suppose $f_s<f_r$. So $r$ is odd.
Let $x:=e(G_n(r,s))=(n-1)\cdot{\bf 1}_r+{\bf 1}_s$.
Consider any $n$-graph $G$ and let $y=e(G)$ be its degree enumerator.
Since $n$ is odd, $G$ must have at least one vertex of even degree $t$,
as otherwise the sum of degrees would have been odd, which is impossible. So $y_t\geq 1$
and $f_t\leq f_s$ by the choice of $s$. Therefore
\begin{eqnarray*}
fy
& = & \sum_{i=0}^{n-1}f_iy_i\ =\ \sum\{f_iy_i:i\neq t\}+f_t(y_t-1)+f_t \\
& \leq & f_r\left(\sum\{y_i:i\neq t\}+(y_t-1)\right)+f_s\ =\ (n-1)f_r+f_s\ =\ fx\ .
\end{eqnarray*}
So $x$ maximizes $f$ over $\E_n$ and since $v$ is the
unique maximizer we have $v=x=(n-1)\cdot{\bf 1}_r+{\bf 1}_s$.

In the other direction, we show that each vector indicated by the theorem is a vertex.
For any even $n$ and any $r$, or any odd $n$ and any even $r$, let $x:=e(G_n(r))=n\cdot{\bf 1}_r$,
and let $f={\bf 1}_r$. Then $fx=n$. Consider any degree enumerator $y=e(G)$ of any $n$-graph $G$
with $y\neq x$. Then $y_r<n$ and hence $fy=y_r<n=fx$.
So $x$ is the unique maximizer of $fz$ over $\E_n$.

Next, for any odd $n$, any odd $r$, and any even $s$,
let $x:=e(G_n(r,s))=(n-1)\cdot{\bf 1}_r+{\bf 1}_s$,
and let $f=2\cdot{\bf 1}_r+{\bf 1}_s$. Then $fx=2(n-1)+1=2n-1$.
Consider any degree enumerator $y=e(G)$ of any $n$-graph $G$ with $y\neq x$.
Since $n$ is odd $G$ must have some vertex of even degree and therefore $y_r\leq n-1$.
If $y_r=n-1$ then $y\neq x$ implies $y_s=0$ and hence
$$fy\ =\ 2y_r+y_s\ =\ 2(n-1)\ <\ 2n-1\ =\ fx\ .$$
If $y_r<n-1$ then
$$fy\ =\ 2y_r+y_s\ \leq\ 2y_r+(n-y_r)\ =\ n+y_r\ <\ 2n-1\ =\ fx\ .$$
So again $x$ is the unique maximizer of $fz$ over $\E_n$.

Thus, the number of extremal degree enumerators is the number of vectors indicated in the theorem,
and it is easy to see that it is indeed $n$ for even $n$ and $({n+1\over2})^2$ for odd $n$.
\epr

\bc{complete_optimization}
Degree sequence optimization over $K_n$ and any $f$ can be done in time $O(n^2)$.
\ec
\bpr
Compute the product $fe$ for each of the $O(n^2)$ vertices $e_n(r)$ and $e_n(r,s)$ of $\E_n$ in
Theorem \ref{complete}, pick $e$ attaining maximum, and output the corresponding $G_n(r)$ or $G_n(r,s)$.

Or faster, along the proof of Theorem \ref{complete}, pick $0\leq r\leq n-1$ with
$f_r=\max\{f_0,\dots,f_{n-1}\}$ and pick $0\leq s\leq n-1$ even with
$f_s=\max\{f_i:0\leq i\leq n-1,\ i\ \mbox{even}\}$. If $n$ is even then the graph $G_n(r)$ is an
optimal solution. If $n$ is odd and $f_s=f_r$ then the graph $G_n(s)$ is an optimal solution.
If $n$ is odd and $f_s<f_r$ then $r$ is odd and the graph $G_n(r,s)$ is optimal.
\epr

We conclude this section with a simple construction showing that the total number of degree enumerators of $n$-graphs
grows exponentially with $n$. In contrast, by Theorem \ref{complete}, of these, the number
of extremal degree enumerators is bounded by a quadratic function of $n$.

\bp{exponential}
The number of degree enumerators of graphs on $[2k+1]$ is at least $2k\choose k$.
\ep
\bpr
Partition the vertex set $V=\{0,\dots,2k\}$ into a $k$-subset $L$ and a $(k+1)$-subset $R$,
say $L$ consists of the odd $i\in V$ and $R$ of the even $i\in V$. For each of the $2k\choose k$
choices of $k$ integers $1\leq s_1<\cdots<s_k\leq 2k$ we construct a graph $G$ on $V$, such
that the degree enumerators $e(G)$ of graphs corresponding to different choices are different.
Let $s_0:=0$ and $s_{k+1}:=2k+1$. Given a choice let $e_i:=s_{i+1}-s_i-1\geq 0$ for $i=0,\dots,k$.
Note that $\sum_{i=0}^k e_i=\sum_{i=0}^k(s_{i+1}-s_i-1)=k$ and pick some partition
$L=\uplus_{i=0}^k L_i$ with $|L_i|=e_i$ for $i=0,\dots,k$. Connect in $G$ each vertex in $L_i$
to some $i$ vertices in $R$. Now add in $G$ a clique on $R$. Then for $i=0,\dots,k-1$ each vertex
in $L_i$ has degree $i$ and all other vertices, in $L_k\uplus R$, have degree at least $k$.
So the prefix of the degree enumerator $e(G)$ consisting of the first $k$ entries satisfies
$(e_0(G),\dots,e_{k-1}(G))=(e_0,\dots,e_{k-1})$. Since $e_0,\dots,e_{k-1}$ determine $s_1,\dots,s_k$
by $s_{i+1}=e_i+s_i+1$ for $i=0,\dots,k-1$, different choices of the $s_i$ give graphs $G$
with different prefixes of $e(G)$ and hence different degree enumerators, completing the proof.
\epr

\section{On extremal bi-enumerators of complete bipartite graphs}

We now consider the complete bipartite graph $K_{m,n}=(U,V,E)$ with $U=\{u_1,\dots,u_m\}$
the set of left vertices, $V=\{v_1,\dots,v_n\}$ the set of right vertices,
and $E=U\times V$ the set of edges. We characterize, for $m=1,2$, and any $n$,
the vertices of the bi-enumerator polytope
$$\B_{m,n}\ =\ \B_{K_{m,n}}\ \subset\ \R^{n+1}\oplus\R^{m+1}\ .$$

\bp{m=1}
For $m=1$, all bi-enumerators are vertices of $\B_{1,n}$, and they are precisely
$${\bf 1}_k\oplus(n-k,k),\ \ k=0,\dots,n\ .$$
Moreover, the degree bi-enumerator polytope $\B_{1,n}$ is an $n$-simplex in
$\R^{n+1}\oplus\R^{1+1}\cong \R^{n+3}$.
\epr
\ep
\bpr
Let $G\subseteq K_{1,n}$ be any subgraph and let $k=d_{u_1}(G)\in[n+1]$ be the degree of
vertex $u_1$. Then $k$ vertices in $V$ have degree $1$ and $n-k$ have degree $0$ so
$a(G)\oplus c(G)={\bf 1}_k\oplus(n-k,k)$. And clearly for any $k\in[n+1]$,
any graph where $u_1$ is connected to precisely $k$ vertices in $V$ has
$d_{u_1}(G)=k$. So the $n+1$ vectors listed above are precisely all
bi-enumerators of subgraphs of $K_{1,n}$. Since each involves a different
unit vector on the left, they are linearly independent and hence also affinely independent.
So $\B_{1,n}$ is an $n$-dimensional simplex in $\R^{n+1}\oplus\R^{1+1}$.
\epr

We proceed with the case of $m=2$ which turns out to be much more complicated. First we
determine all bi-enumerators of all subgraphs $G\subseteq K_{2,n}$.
Consider any subgraph $G\subseteq K_{2,n}$. Let $i:=d(u_1)$ and $j:=d(u_2)$ and
without loss of generality assume $0\leq i\leq j\leq n$. Let $N(u_1),N(u_2)\subseteq V$ be the sets
of neighbors of $u_1,u_2$ in $G$ and let $k:=|N(u_1)\cap N(u_2)|$ be the size of their
intersection. So $i=|N(u_1)|\leq|N(u_2)|=j$ and $\max\{0,(i+j)-n\}\leq k\leq i$.
Then in $V$ there are $|[n+1]\setminus(N(u_1)\cup N(u_2))|= n-(i+j)+k$ vertices of degree $0$,
$|(N(u_1)\setminus N(u_2))\uplus(N(u_2)\setminus N(u_1))|=(i+j)-2k$ vertices of degree $1$, and
$|N(u_1)\cap N(u_2)|=k$ vertices of degree $2$. Therefore the degree bi-enumerator of $G$ satisfies
$$b(G)\ =\ a(G)\oplus c(G)\ =\ b_n(i,j,k)\ :=\ ({\bf 1}_i+{\bf 1}_j)\oplus(n-(i+j)+k,(i+j)-2k,k)\ .$$
And, we can construct a graph $G_n(i,j,k)\subseteq K_{2,n}$ with $b(G_n(i,j,k))=b_n(i,j,k)$ for any
$$i=0,1,\dots,n,\quad j=i,\dots,n,\quad k=\max\{0,(i+j)-n\},\dots,i$$
by picking some $i$-subset $N_1\subseteq V$ and some $j$-subset $N_2\subseteq V$ such that
$|N_1\cap N_2|=k$, and connecting $u_r$ to all vertices in $N_r$ for $r=1,2$.
We therefore have the following proposition.

\bp{m=2_bi-enumerators}
For $m=2$, the bi-enumerators of all $G\subseteq K_{2,n}$ are precisely the following,
$$b_n(i,j,k)\ =\ b(G_n(i,j,k))\ =\ ({\bf 1}_i+{\bf 1}_j)\oplus(n-(i+j)+k,(i+j)-2k,k)\ ,$$
$$i=0,1,\dots,n,\quad j=i,\dots,n,\quad k=\max\{0,(i+j)-n\},\dots,i\ .$$
Their number is ${(n+2)(n+4)(2n+3)\over 24}$ for even $n$
and ${(n+1)(n+3)(2n+7)\over 24}$ for odd $n$, so is cubic in $n$.
\ep

\bpr
The first statement follows by the discussion and construction of the $G_n(i,j,k)$ above.

Next we count the number of bi-enumerators which we denote by $a(n)$. For even $n$ it is
\begin{eqnarray*}
a(n)&=&|\{(i,j,k):i=0,1,\dots,n,\ j=i,\dots,n,\ k=\max\{0,(i+j)-n\},\dots,i\}| \\
&=&|\{(i,j,k):i=0,\dots,{n\over2},\ j=i,\dots,n-i,\ k=0,\dots,i\}| \\
&+&|\{(i,j,k):j={n\over2}+1,\dots,n,\ i=n-j+1,\dots,j,\ k=(i+j)-n,\dots,i\}| \\
&=&\sum_{i=0}^{n\over2}(n-2i+1)(i+1)+\sum_{j={n\over2}+1}^n(2j-n)(n-j+1)
\ =\ {(n+2)(n+4)(2n+3)\over 24}\ ,
\end{eqnarray*}
as claimed, where the second equality follows by splitting the set of triples to those
with $i+j\leq n$ and those with $i+j>n$, and the last equality is obtained by using
the well known formulas for $\sum_{r=1}^s r$ and $\sum_{r=1}^s r^2$. For odd $n$ a similar
computation applies and omitted.
\epr

We now show that the extremal bi-enumerators of $K_{2,n}$ are the bi-enumerators of those graphs
$G_n(i,j,k)$ where either both vertices on the left have the same degree $i=j$, and their
neighborhoods have either minimum or maximum possible intersection $k=\max\{0,2i-n\}$ or $k=i$,
or one vertex on the left has degree $i\leq\lfloor{n-1\over2}\rfloor$,
the other vertex on the left has degree $j\geq\lceil{n+1\over2}\rceil$,
and their neighborhoods have minimum possible intersection $k=\max\{0,2i-n\}$.

\bt{m=2_vertices}
The vertices of the degree bi-enumerator polytope $\B_{2,n}$ are the following:
\begin{enumerate}
\item
$b_n(i,i,k)=2\cdot{\bf 1}_i\oplus(n-2i+k,2i-2k,k)$
for $i=0,\dots,n$ and $k=\max\{0,2i-n\},i$;
\item
$b_n(i,j,k)$ for $i=0,\dots,\lfloor{n-1\over2}\rfloor$ and $j=\lceil{n+1\over2}\rceil,\dots,n$
where $k=\max\{0,(i+j)-n\}$.
\end{enumerate}
The number of vertices is $({n\over2})^2+2n$ for even $n$
and $({n+1\over2})^2+2n$ for odd $n$, so is quadratic in $n$.\break
\et
\bpr
We need to determine which of the bi-enumerators given in Proposition \ref{m=2_bi-enumerators}
are vertices. Since $n$ is understood, for brevity we denote throughout the proof
$b(i,j,k):=b_n(i,j,k)$. Consider any $0\leq i\leq j\leq n$. There are bi-enumerators $b(i,j,k)$
for $\max\{0,(i+j)-n\}\leq k\leq i$. If there is any $k$ with $\max\{0,(i+j)-n\}<k<i$ then
$b(i,j,k)={1\over2}(b(i,j,k-1)+b(i,j,k+1))$, so $b(i,j,k)$ is not a vertex.
So among the bi-enumerators $b(i,j,k)$ corresponding to the pair $i,j$, at most two
could potentially be vertices, those with $k=\max\{0,(i+j)-n\}$ and $k=i$.

Consider first the case $i=j$. Let $f:={\bf 1}_i\in\R^{n+1}$ and $g:=0\in\R^{2+1}$.
Then the face of $\B_{2,n}$ on which the inner product $\langle f\oplus g,x\oplus y\rangle$
is maximized is precisely the convex hull of the bi-enumerators $b(i,i,k)$ for all relevant $k$.
By what said above, we conclude that $b(i,i,\max\{0,2i-n\})$ and $b(i,i,i)$ are the vertices of that
face and hence of $\B_{2,n}$. Note that for $i=0$ or $i=n$ these collapse to a single vertex
$b(i,i,i)$, for $0<i\leq {n\over 2}$ these are the two distinct vertices
$b(i,i,0)$ and $b(i,i,i)$, and for ${n\over 2}\leq i<n$ these are the two distinct vertices
$b(i,i,2i-n)$ and $b(i,i,i)$. So the total number of vertices of the form $b(i,i,k)$
is $1+(n-1)\cdot 2+1=2n$.

Now consider the case $i<j$. Consider three subcases. First, suppose $0\leq i<j\leq{n\over2}$.
Then both of the potential vertices $b(i,j,0)$ and $b(i,j,i)$ are not vertices, since we have
$$b(i,j,0)\ =\ {1\over2}(b(i,i,0)+b(j,j,0))\ ,\quad
b(i,j,i)\ =\ {1\over2}\cdot b(i,i,i)+{j-i\over 2j}\cdot b(j,j,0)+{i\over 2j}\cdot b(j,j,j)\ .$$
For instance, for $n=5$, $i=1$, $j=2$, we have
$b(1,2,1)={1\over2}b(1,1,1)+{1\over4}b(2,2,0)+{1\over4}b(2,2,2)$,
$$011000\oplus311\ =\ {1\over2}(020000\oplus401)+
{1\over4}(002000\oplus140)+{1\over4}(002000\oplus302)\ .$$
Second, suppose ${n\over2}\leq i<j\leq n$. Then $b(i,j,(i+j)-n)$ and $b(i,j,i)$
are not vertices, since
\begin{eqnarray*}
b(i,j,(i+j)-n)&=&{1\over2}(b(i,i,2i-n)+b(j,j,2j-n))\ ,\\
b(i,j,i)&=&{n-j\over2(n-i)}\cdot b(i,i,i)+
{j-i\over 2(n-i)}\cdot b(i,i,2i-n)+{1\over 2}\cdot b(j,j,j)\ .
\end{eqnarray*}
Third, suppose $0\leq i\leq\lfloor{n-1\over2}\rfloor$ and $\lceil{n+1\over2}\rceil\leq j\leq n$.
For each such pair $i,j$ there are again potentially at most two vertices $b(i,j,k)$ with
$k=\max\{0,(i+j)-n\}$ and $k=i$. We show that for each such pair $i,j$ we have that $b(i,j,k)$
is indeed a vertex for $k=\max\{0,(i+j)-n\}$, but is not a vertex for $k=i$,
except for $i=0$ or $j=n$, in which case $i=\max\{0,(i+j)-n\}$.

We begin by showing that $b(i,j,i)$ is not a vertex for
$0<i\leq\lfloor{n-1\over2}\rfloor$ and $\lceil{n+1\over2}\rceil\leq j<n$.
First suppose $j\leq 2i$. Then $b(i,j,i)$ is not a vertex
as it is the convex combination
$$b(i,j,i)\ =\ {j-i\over 2i}\cdot b(i,i,0)+{2i-j\over 2i}\cdot b(i,i,i)+{1\over2}\cdot b(j,j,j)\ .$$
Next suppose $j\geq 2i$ and $i+j\leq n$. Then $b(i,j,i)$ is not a vertex
as it is the convex combination
$$b(i,j,i)\ =\ {i\over j}\cdot b(i,i,0)+{i\over j}\cdot b(j,j,j)+{j-2i\over j}\cdot b(i,j,0)\ .$$
For instance, for $n=5$, $i=1$, $j=3$, we have
$b(1,3,1)={1\over3}b(1,1,0)+{1\over3}b(3,3,3)+{1\over3}b(1,3,0)$,
$$010100\oplus221\ =\ {1\over3}(020000\oplus320)+
{1\over3}(000200\oplus203)+{1\over3}(010100\oplus140)\ .$$
Last suppose $j\geq 2i$ and $i+j\geq n$. Then $b(i,j,i)$ is not a vertex
as it is the convex combination
$$b(i,j,i)\ =\ {n-j\over n-i}\cdot b(i,i,i)+{n-j\over n-i}\cdot b(j,j,2j-n)
+{((i+j)-n)+(j-2i)\over n-i}\cdot b(i,j,(i+j)-n)\ .$$

We proceed to show that $b(i,j,k)$ is indeed a vertex for $0\leq i\leq\lfloor{n-1\over2}\rfloor$,
$\lceil{n+1\over2}\rceil\leq j\leq n$, and $k=\max\{0,(i+j)-n\}$.
Support first $i+j\leq n$. Let $f:=N({\bf 1}_i+{\bf 1}_j)\in\R^{n+1}$ where $N$ is a
suitably large positive integer and $g:=(0,2j-n,(i+j)-n)\in\R^{2+1}$.
We claim that the inner product $\langle f\oplus g,b(i',j',k')\rangle$ is
maximized over $\B_{2,n}$ uniquely at $b(i,j,0)$. Indeed, we have
\begin{eqnarray*}
&\langle f\oplus g,b(i',j',k')\rangle&
=\ \langle N({\bf 1}_i+{\bf 1}_j),{\bf 1}_{i'}+{\bf 1}_{j'} \rangle \\
&&+\ \langle (0,2j-n,(i+j)-n),(n-(i'+j')+k',(i'+j')-2k',k') \rangle \\
&&\hskip-0.65cm\left\{
  \begin{array}{ll}
    \ =\ 2N+(2j-n)(i+j), & \hbox{if $(i',j',k')=(i,j,0)$;} \\
    \ \leq\ N+3n^2, & \hbox{if $(i',j')\neq(i,j)$ and $i'<j'$;} \\
    \ \leq\ 3n^2, & \hbox{if $i'=j'\neq i,j$;} \\
    \ =\ 2N+(2j-n)2i, & \hbox{if $i'=j'=i, k'=0$;} \\
    \ =\ 2N+((i+j)-n)i, & \hbox{if $i'=j'=i, k'=i$;} \\
    \ =\ 2N+(2j-n)(i+(n-j)), & \hbox{if $i'=j'=j, k'=2j-n$;} \\
    \ =\ 2N+((i+j)-n)j, & \hbox{if $i'=j'=j, k'=j$.}
  \end{array}
\right.
\end{eqnarray*}
The claim now follows since $i+j\leq n$, $2i<i+j$, $n-j<j$,
and since $N$ is sufficiently large.

Now assume $i+j>n$. Let $f:=N({\bf 1}_i+{\bf 1}_j)\in\R^{n+1}$
where $N$ is a suitably large positive integer and $g:=(0,2((i+j)-n),2((i+j)-n)-1)\in\R^{2+1}$.
We now claim that the product $\langle f\oplus g,b(i',j',k')\rangle$ is maximized
over $\B_{2,n}$ uniquely at $b(i,j,(i+j)-n)$. Indeed, we now have
\begin{eqnarray*}
&\langle f\oplus g,b(i',j',k')\rangle&
=\ \langle N({\bf 1}_i+{\bf 1}_j),{\bf 1}_{i'}+{\bf 1}_{j'}\rangle \\
&&+\ 2((i+j)-n)((i'+j')-2k')+(2((i+j)-n)-1))k' \\
&&=\ \langle N({\bf 1}_i+{\bf 1}_j),{\bf 1}_{i'}+{\bf 1}_{j'}\rangle +\ 2((i+j)-n)((i'+j')-k')-k' \\
&&\hskip-0.65cm\left\{
  \begin{array}{ll}
    \ =\ 2N+((i+j)-n)(2n-1), & \hbox{if $(i',j',k')=(i,j,(i+j)-n)$;} \\
    \ \leq\ N+2n^2, & \hbox{if $(i',j')\neq(i,j)$ and $i'<j'$;} \\
    \ \leq\ 2n^2, & \hbox{if $i'=j'\neq i,j$;} \\
    \ =\ 2N+((i+j)-n)4i, & \hbox{if $i'=j'=i, k'=0$;} \\
    \ =\ 2N+((i+j)-n)2i-i, & \hbox{if $i'=j'=i, k'=i$;} \\
    \ =\ 2N+((i+j)-n)2n-(2j-n), & \hbox{if $i'=j'=j, k'=2j-n$;} \\
    \ =\ 2N+((i+j)-n)2j-j, & \hbox{if $i'=j'=j, k'=j$.}
  \end{array}
\right.
\end{eqnarray*}
The claim again follows since $i<j$, $4i\leq2(n-1)$,
and since $N$ is sufficiently large.

For instance, for $n=5$, $i=2$, $j=4$, $f:=N({\bf 1}_2+{\bf 1}_4)\in\R^{5+1}$
with suitably large $N$, and $g:=(0,2,1)$, the product
$\langle f\oplus g,b(i',j',k')\rangle$ is maximized over $\B_{2,5}$ uniquely at $b(2,4,1)$.
\epr

\bc{bipartite_optimization}
Degree bi-optimization over $K_{2,n}$ and any $f$ and $g$ can be done in time $O(n^2)$.\break
\ec
\bpr
Compute the product $\langle f\oplus g,b_n(i,j,k)\rangle$ for each of the
$O(n^2)$ vertices of $\B_{2,n}$ appearing in Theorem \ref{m=2_vertices},
pick one attaining maximum, and output the corresponding $G_n(i,j,k)$.
\epr

\section*{Acknowledgments}

The author thanks Frédéric Meunier for discussions on this work and in particular on the hardness
of the degree sequence optimization problem. The author was partially supported by the
National Science Foundation under Grant No. DMS-1929284 while he was in residence at
the Institute for Computational and Experimental Research in Mathematics, Brown University,
during the semester program {\em Discrete Optimization: Mathematics, Algorithms, and Computation}.
The author was also partially supported by the Dresner chair at the Technion.

\end{document}